\documentclass[titlepage,11pt]{article}
\usepackage{latexsym}
\usepackage{amsfonts}
\usepackage{amsmath}
\newcommand\blackslug{\hbox{\hskip 1pt \vrule width 4pt height 8pt depth 1.5pt
        \hskip 1pt}}
\newcommand\bbox{\hfill \quad \blackslug \bigbreak}

\def\l{,\ldots,}

\title{Immersion in four-edge-connected graphs}
\author{Maria Chudnovsky\thanks{Supported by NSF grants IIS-1117631 and DMS-1001091.}\\
Columbia University, New York, NY 10027, USA
\\
\\
Zden\v{e}k Dvo\v{r}\'ak\thanks{Supported by the Center of Excellence of the Institute for
Theoretical  Computer Science, Prague,
project P202/12/G061 of Czech Science Foundation.}\\
Charles University, Prague, Czech Republic
\\
\\
Tereza Klimo\v{s}ov\'a\thanks{This author was a student at Charles University
until September 2012. The work leading to this invention has received funding from the European
Research Council  under the European Union's Seventh Framework Programme 
(FP7/2007-2013)/ERC  grant agreement no.~259385.}\\
University of Warwick, Coventry CV4 7AL, UK
\\
\\
Paul Seymour\thanks{Supported by ONR grant N00014-10-1-0680 and NSF grant DMS-0901075.}\\
Princeton University, Princeton, NJ 08544, USA}

\date{January 3, 2013; revised \today}

\newtheorem{thm}{}[section]

\newcommand{\Proof}{\noindent{\bf Proof.}\ \ }

\begin{document}
\maketitle
\begin{abstract}
Fix $g>1$. Every graph of large enough tree-width contains a $g\times g$ grid as a minor; but here we prove that every four-edge-connected
graph of large enough tree-width contains a $g\times g$ grid as an immersion (and hence contains any fixed graph with maximum degree at most four as an immersion). 
This result has a number of applications.
\end{abstract}

\section{Introduction}
Let $G,H$ be graphs. (All graphs in this paper are finite, possibly with loops or parallel edges.) 
An {\em immersion} of $H$ in $G$ is a map $\eta$, with domain $V(H)\cup E(H)$, mapping each vertex of 
$H$ to a vertex of $G$, and each edge of $H$ to a path or cycle of $G$, satisfying the following:
\begin{itemize}
\item $\eta(u)\ne \eta(v)$ for all distinct $u,v\in V(H)$
\item for each $e\in E(H)$ with distinct ends $u$ and $v$, $\eta(e)$ is a path of $G$ with ends $\eta(u),\eta(v)$
\item for each loop in $H$ with end $v$, $\eta(e)$ is a cycle of $G$ passing through $\eta(v)$
\item for $v\in V(H)$ and $e\in E(H)$, if $e$ is not incident with $v$ in $H$ then $\eta(v)\notin V(\eta(e))$
\item for all distinct $e,f\in E(H)$, $E(\eta(e)\cap \eta(f)) = \emptyset$.
\end{itemize}
(Note in particular the fourth condition above; this relation is normally called ``strong immersion'', but we omit ``strong'' since we do not need weak immersion in this paper.)
If there is an immersion of $H$ in $G$, we say that ``$H$ can be immersed in $G$'' and ``$G$ contains $H$ as an immersion'' (or just ``$G$ immerses $H$'').
If in addition, for all distinct $e,f\in E(H)$, every vertex of $\eta(e)\cap \eta(f)$ is equal to $\eta(v)$ for some $v\in V(H)$ incident in $H$ 
with both $e$ and $f$, then $\eta$ is called a {\em subdivision map} of $H$ in $G$.

If $g>1$ is an integer, the $g\times g$ {\em grid} is a graph with vertex set $\{v_{ij}\;:1\le i,j\le g\}$, where $v_{ij}$ is adjacent to $v_{i'j'}$ 
if $|i-i'|+|j-j'| = 1$. We denote this graph by $J_g$.
Evidently any graph of maximum degree at most four can be immersed in $J_g$, if $g$ is sufficiently large (map the vertices far apart, and then route the edges
so that no three pass through the same vertex, and any two that share a vertex cross there). 

A {\em tree-decomposition} of a graph $G$ is a pair $(T,(W_t\;:t\in V(T)))$, such that 
\begin{itemize}
\item $T$ is a tree
\item  $W_t\subseteq V(G)$ for each $t\in V(T)$
\item $V(G) = \cup (W_t\;:t\in V(T))$
\item for every edge $uv$ of $G$, there exists $t\in V(T)$ with $u,v\in W_t$
\item for $t,t',t''\in V(T)$, if $t'$ belongs to the path of $T$ between $t$ and $t''$, then $W_t\cap W_{t''}\subseteq W_{t'}$.
\end{itemize}
We call $\max(|W_t|-1\;:t\in V(T))$ the {\em width} of the tree-decomposition, and say that 
$G$ has {\em tree-width} $k$ if $k$ is minimum such that $G$ admits a tree-decomposition of width $k$.

We say that $H$ is a {\em minor} of $G$ if a graph isomorphic
to $H$ can be obtained from a subgraph of $G$ by contracting edges. The following is well-known~\cite{GM5}:
\begin{thm}\label{gridminor}
For all $g>1$ there exists $k$ such that every graph with tree-width at least $k$ contains $J_g$ as a minor.
\end{thm}
(Note that this is sharp in the sense that for all $k$ there exists $g$ such that no graph of tree-width less than $k$ contains $J_g$ as a minor.) 
In this paper we prove a similar result for immersion, the following. (Two versions of this result were found independently by two subsets of the authors,
and one of these versions appears in \cite{klimosova}.)

\begin{thm}\label{main1}
For all $g>1$ there exists $k\ge 0$ such that every four-edge-connected graph with tree-width at least $k$ contains $J_g$ as an immersion.
\end{thm}

This is not exactly an analogue of \ref{gridminor}, because it is not sharp in the same sense. It is not true that for all $k$ there exists $g$ such that no 
four-edge-connected graph of tree-width less than $k$ contains $J_g$ as a minor. To see this, let $G$ be obtained from a star with $g^2+1$ vertices
by replacing each edge by four parallel edges. Then $G$ has tree-width one, and yet immerses $J_g$. Thus large tree-width is too strong a requirement.
One might think that we should measure the width by decomposing the graph with a tree-structure of edge-cutsets of bounded size, instead of vertex-cutsets (which is
essentially what tree-width does); but this is now too {\em weak}; a two-vertex graph with many parallel edges has large width under this measure, and yet does
not immerse $J_2$. So the correct concept is somewhere between the two, and getting it right is beyond the scope of this paper (see
\cite{wollan}). (Note that these problems 
arise since there is no bound on the maximum degree of $G$; if we bound the maximum degree, then the vertex-cut and edge-cut versions of tree-width become equivalent,
and \ref{main1} becomes sharp in the desired sense.)

Let $h\ge 2$ be even. An {\em elementary wall of height $h$} is a graph whose vertex set can be labeled
$$\{v_{ij}\;:1\le i\le h+1,1\le j\le 2h+2, (i,j)\ne (1,2h+2), (h+1,1)\}$$
where distinct vertices $v_{ij}, v_{i'j'}$ are adjacent if either
\begin{itemize}
\item $i = i'$ and $|j'-j| = 1$, or
\item $j = j'$ and $|i'-i| = 1$ and $\min(i,i')+j$ is even.
\end{itemize}
A {\em wall of height $h$} is a subdivision of an elementary wall of height $h$, and a {\em wall contained in a graph $G$} (or just a ``wall in $G$'') 
means a subgraph of $G$ that is a wall.

The advantage of walls is that they permit us to state a version of \ref{gridminor} using subgraphs instead of minors, the following (this is easy to see):
\begin{thm}\label{wallminor}
For all even $h\ge 2$ there exists $k$ such that every graph with tree-width at least $k$ contains a wall of height $h$.
\end{thm}

Thus the following is equivalent to \ref{main1}:
\begin{thm}\label{main2}
For all $g>1$ there exists $h\ge 2$, even, such that every four-edge-connected graph containing a wall of height $h$ contains $J_g$ as an immersion.
\end{thm}

This admits two strengthenings, which we need to describe next. First, 
with the usual labelling of the vertex set of a wall, the vertices $v_{i,2i}\;(2\le i\le h)$
are called its {\em diagonal} vertices. We will be able to replace the ``four-edge-connected'' hypothesis with a weaker hypothesis that a large number of the diagonal 
vertices are pairwise four-edge-connected. Second
we want to show that if we start with a large wall, we get a large grid immersion that is in some sense
``close'' to the wall. More precisely, we insist that the immersion map each vertex of the grid to one of the diagonal vertices in our $4$-edge-connected set.
If $S\subseteq V(G)$, and $\eta$ is an immersion of $H$ in $G$, we say $\eta$ is {\em $S$-rooted} if $\eta(v)\in S$ for each $v\in V(H)$. 
The following version of \ref{main2} incorporates both the strengthenings just discussed, and is the main result of the paper.
\begin{thm}\label{main3}
For all $g>1$ there exists $b\ge 0$, with the following property. Let $W$ be a wall contained in a graph $G$, 
and let $S$ be a set of diagonal vertices of $W$, pairwise $4$-edge-connected in $G$, and with $|S|\ge b$.
Then there is an $S$-rooted immersion of $J_g$ in $G$.
\end{thm}

\section{Applying lemmas from graph minors}

Let $W$ be a wall of height $h$, and let $W_0$ be the elementary wall of which $W$ is a subdivision. Label the vertices of $W_0$ as in the definition.
Each edge of $W_0$
corresponds to a path of $W$, and we call such a path a {\em branch} of $W$. 
Choose $h'$ even with $2\le h'\le h$, and choose four integers $i_1,i_2,j_1,j_2$, satisfying
\begin{itemize}
\item $i_1,j_1$ are odd 
\item $1\le i_1\le i_2\le h+1$ and $1\le j_1\le j_2\le 2h$
\item $i_2-i_1 = h'$ and $j_2-j_1 = 2h'-1$. 
\end{itemize}
Then the subgraph of $W_0$ induced on the vertex set 
$$\{v_{i,j}: i_1\le i\le i_2, j_1\le j\le j_2, (i,j)\ne (i_1,j_2), (i_2,j_1)\}$$ 
is an elementary wall of height $h'$, and we call such a wall an {\em elementary subwall} of $W_0$. The corresponding subgraph of $W$
is called a {\em subwall} of $W$.

Let $W$ be a wall, in a graph $G$, and let $v$ be a diagonal vertex of $W$. The {\em surround} of $v$ means the set of vertices $u$ of $W$ such that
either $u = v$, or there is a path of $W$ between $u,v$ in which every vertex different from $v$ has degree two in $W$. 
A {\em fin} for $W$ means a triple $(s,F,t)$ such that
\begin{itemize}
\item $s$ is a diagonal vertex of $W$, 
\item $t \in V(W)$ does not belong to the surround of $s$, and
\item $F$ is a path in $G$ with ends $s,t$, edge-disjoint from $W$.
\end{itemize}
We also call this a {\em fin at $s$}. We say that 
$(W,(s_i,F_i,t_i)_{1\le i\le b})$
 is a {\em fin system in $G$ } if
\begin{itemize}
\item $W$ is a wall in $G$
\item $s_1\l s_b$ are diagonal vertices of $W$, all different
\item for $1\le i\le b$, $(s_i,F_i,t_i)$ is a fin for $W$, and
\item for $1\le i,j\le b$, if $i\ne j$ then $s_i\notin V(F_j)$.
\end{itemize}

In this section we prove the following.
\begin{thm}\label{fins}
For all $g>1$ there exists an integer $b\ge 0$, with the following property. Let 
$$(W,(s_i,F_i,t_i)_{1\le i\le b})$$
be a fin system in $G$. 
Then there is an $\{s_1\l s_b\}$-rooted immersion of $J_g$ in $G$.
\end{thm}

The proof is in several steps. Some of the proofs of the steps are merely sketches, because they are standard applications of methods of the graph minors papers (rather
straightforward, since the underlying graph is a wall), and to import all the definitions and theorems of the corresponding graph minors papers and to apply them 
precisely would take a considerable amount
of space (and we suspect would not improve clarity).

Let $W$ be an elementary wall of height $h$; then 
there is a drawing of $W$ in the plane so that all finite regions have boundary of length six (and such a drawing is unique up to homeomorphisms of the plane).
We call this the {\em standard drawing}. A standard drawing of a wall is obtained by subdividing edges of the standard drawing of the corresponding 
elementary wall, and we call the finite regions {\em bricks}.

Let $W$ be a wall with its standard drawing.
If $s,t$ are points of the plane, we define
their {\em distance}
$d(s,t)$ to be zero if $s = t$, and otherwise to be the minimum of the number of points of $F$ in the drawing, taken over all subsets $F$ in the plane homeomorphic
to $[0,1]$ with ends $s$ and $t$. The {\em perimeter} of $W$ is the cycle bounding the infinite region.

We begin with
\begin{thm}\label{longjumps}
Let $g>1$. Then there exist integers $a_1,b_1>0$ with the following property. Let $G$ be a graph 
and  $(W,(s_i,F_i,t_i)_{1\le i\le b_1})$
be a fin system in $G$.
Suppose in addition that:
\begin{itemize}
\item the vertices $s_1\l s_{b_1}, t_1\l t_{b_1}$ are all pairwise at distance at least $a_1$, and
\item the paths $F_1\l F_{b_1}$ are pairwise edge-disjoint.
\end{itemize}
Then there is an $\{s_1\l s_{b_1}\}$-rooted immersion of $J_g$ in $G$.
\end{thm}
\Proof We may assume that $g$ is even (by replacing $g$ by $g+1$ if necessary). Let $a_1$ be large (in terms of $g$). Let $n = g^2$ and $b_1 = n+1$.
Let $J_g$ have vertex set $\{j_1\l j_n\}$ say. Since $g$ is even, there is a perfect matching $M$ in $J_g$. 

Now let $G$ and $(W,(s_i,F_i,t_i)_{1\le i\le b_1})$ be as in the theorem. Since $s_1\l s_{b_1}, t_1\l t_{b_1}$ are pairwise at distance at least $a_1$, 
at most one of them has distance less than $a_1/2$ 
of a vertex of the perimeter of $W$; so we may assume that none of $s_1\l s_{n}, t_1\l t_{n}$ is within distance $(a_1-1)/2$
of the perimeter. 

Now each $s_i$ has degree three in $W$; let the three neighbours of $s_i$ in $W$ be $x_{i,1},x_{i,2}, x_{i,3}$, enumerated in clockwise order around $s_i$, 
with an arbitrary first vertex. 
For $1\le i\le n$, let the edges of $J_g$ incident with $j_i$ be $e_{i,1}\l e_{i,k_i}$ (where $k_i$ is the degree of $j_i$ in $J_g$), enumerated in clockwise order around $j_i$,
such that $e_{i,k_i}\in M$.
If $a_1$ is large enough (in terms of $g$), then by theorem 7.5 of~\cite{GM7} applied to $W\setminus \{s_1\l s_{n}\}$,
(and compare theorem 4.4 of~\cite{GM12} for a similar situation), we deduce that
for each edge $e$ of $J_g$ there is a path $Q_e$ of $W\setminus \{s_1\l s_n\}$ satisfying the following, where $e$ has ends $j_h, j_i$ say in $J_g$:
\begin{itemize}
\item if $e \in M$ then $Q_e$ has ends $t_h,t_i$;
\item if $e\notin M$, let $e = e_{h,p} = e_{i,q}$; then the ends of $Q_e$ are $x_{h,p}$ and $x_{i,q}$;
\item the paths $Q_e\;(e\in E(J_g))$ are pairwise vertex-disjoint.
\end{itemize}
For each $e\in M$ with ends $j_h,j_i$ say, let $\eta(e)$ be a path between  $s_h,s_i$ in the the union of the paths $F_h,Q_e$ and $F_i$.
For each $e\in E(J_g)\setminus M$ with ends $j_h,j_i$ say, let $p,q$ satisfy $e = e_{h,p} = e_{i,q}$, and let $\eta(e)$ be the path between $s_h,s_i$
formed by the union of $s_hx_{h,p}$,
$x_{i,q}s_i$ to $Q_e$.
Let $\eta(j_i) = s_i\;(1\le i\le n)$; then $\eta$ is the desired immersion.
This proves \ref{longjumps}.~\bbox

Next we need:

\begin{thm}\label{outblob}
Let $g>1$, and let $a_1,b_1$ satisfy \ref{longjumps}. Let $G$ be a graph, and let $W$ be a wall in
$G$. Let $X$ be a connected subgraph of $G$, edge-disjoint from $W$. Let $S$ be a set of $2b_1$ diagonal
vertices of $W$, pairwise at distance at least $a_1$, each in $V(X)$ and with degree one in $X$.
Then there is an $S$-rooted immersion of $J_g$ in $G$.
\end{thm}
\Proof
Let $G,W,X,S$ be as in the theorem. We may assume that $|S| = 2b_1$.
Since $X$ is connected, there is a spanning tree $T$ of $X$.
Since $|S|$
is even, it is possible to pair up the vertices in $S$ such that there are pairwise edge-disjoint paths of $T$ joining the pairs. 
This provides a fin system satisfying the hypotheses of \ref{longjumps}, and the result follows.~\bbox

\begin{thm}\label{inblob}
Let $g>1$. Then there exists integers $a_2,b_2$ with the following property. Let $G$ be a graph, and let 
$(W,(s_i,F_i,t_i)_{1\le i\le b_2})$
be a fin system in $G$.
Suppose that 
\begin{itemize}
\item $s_1\l s_{b_2}$ are pairwise at distance at least $a_2$, 
\item the vertices $s_i,t_i$ are at distance at least $a_2$ for $1\le i\le b_2$, and
\item the paths $F_1\l F_{b_2}$ are pairwise edge-disjoint.
\end{itemize}
Then there is an $\{s_1\l s_{b_2}\}$-rooted immersion of $J_g$ in $G$.
\end{thm}
\Proof
Let $a_2,b_2$ satisfy $a_2\ge 2a_1$ and $b_2 \ge 6(2b_1+1)b_1$, where $a_1,b_1$ satisfy \ref{longjumps}. Let 
$G$ and 
$$(W,(s_i,F_i,t_i)_{1\le i\le b_2})$$ 
be as in the theorem. 

Since for each $j$ there is at most one $s_i$ with distance less than $a_2/2$ to $t_j$, we may assume that $d(s_i,t_j)\ge a_2/2$ for all $i,j\le n_1$,
where $n_1\ge b_2/3\ge 2(2b_1+2)b_1$.
Suppose that at least $b_1$ of the $t_i$'s pairwise have distance
at least $a_1$; then the result follows from \ref{longjumps}. So for some $j\le n_1$, there are at least $n_1/b_1$ values 
of $i\le n_1$ such that $d(t_i,t_j)<a_1$. Let $n_2 =  2(2b_1+2)$; then $n_1/b_1\ge n_2$, 
and we may assume that $d(t_i,t_1)<a_1$ for $1\le i\le n_2$.

Now $W$ is a wall, and $s_1\l s_{n_2}$ are diagonal vertices of it. There are therefore two subwalls of $W$, say $W_1, W_2$, pairwise disjoint, such that
for $i = 1,2$, exactly $n_2/2$ of $s_1\l s_{n_2}$ are diagonal vertices of $W_i$. 
From the symmetry we may assume that $t_1\notin V(W_1)$. Since  $d(t_i,t_1)<a_1$ for $1\le i\le n_2$, and the $s_i$'s pairwise
have distance at least $a_2$, there is a subwall $W'$ of $W_1$ and hence of $W$ (obtained from $W_1$ by removing an appropriate border) 
such that at least $n_2/2-2\ge 2b_1$ of $s_1\l s_{n_2}$ are diagonal vertices of it, and such that none of 
the corresponding $t_i$'s belong to $W'$. But $W\setminus V(W')$ is connected, and 
the result follows from \ref{outblob}, taking $X$ to be the union of $W\setminus V(W')$ and all the paths $F_i$ with $t_i\in V(W)\setminus V(W')$.~\bbox

\pagebreak

\begin{thm}\label{shortjumps}
Let $g>1,c\ge 0$. Then there exist integers $a_3,b_3\ge 0$ with the following property. Let $G$ be a graph, and let 
$(W,(s_i,F_i,t_i)_{1\le i\le b_3})$
be a fin system in $G$.
Suppose that 
\begin{itemize}
\item $s_1\l s_{b_3}$ are pairwise at distance at least $a_3$
\item $s_i,t_i$ are at distance at most $c$ for $1\le i\le b_3$, and
\item $F_1\l F_{b_3}$ are pairwise edge-disjoint.
\end{itemize}
Then there is an $\{s_1\l s_{b_3}\}$-rooted immersion of $J_g$ in $G$.
\end{thm}
\Proof
Let $n = g^2$.
Let $J_g$ have vertex set $\{j_1\l j_n\}$ say. 
For $1\le i\le n$, let the edges of $J_g$ incident with $j_i$ be $e_{i,1}\l e_{i,k_i}$ (where $k_i$ is the degree of $j_i$ in $J_g$), enumerated in clockwise order around $j_i$.

Let $a_3,b_3$ be big (in terms of $g,c$). 
There is at most one $i$ such that $s_i$ has distance less than $a_3/2$ to some vertex of the perimeter, so we may assume that
$s_i$ has distance at least $a_3/2$ to the perimeter, for $1\le i\le b_3-1$.

Let $1\le i\le b_3-1$. Then there is a subwall $W_i$ of $W$ of height $c$ containing both $s_i,t_i$. Since we choose
$a_3$ much greater than $c$, these subwalls are pairwise disjoint. Let $W'$ be obtained from $W$ by deleting all vertices of $W_i$ for $1\le i\le b_3-1$, and all internal vertices
of branches of $W$ with an end in $W_i$. Let $R_i$ be the region of $W'$ in which $W_i$ was drawn. For $1\le i\le b_3-1$, choose four (distinct) vertices of the boundary of 
$R_i$, say  $x_{i,1}\l x_{i,4}$ in clockwise order around $R_i$, such that there are four paths $B(i,1)\l B(i,4)$ of $W\cup F_i$ from $s_i$ to $x_{i,1}\l x_{i,4}$ respectively,
pairwise edge-disjoint, each with only its final vertex in the boundary of $R_i$.

For $n+1\le i\le b_4$, let us add the edges $x_{i,1}x_{i,3}$ and $x_{i,2}x_{i,4}$ to $W'$,
forming $W''$ say. By theorem 4.5 of~\cite{GM12}, if $b_4-n$ is sufficiently large (in terms of $n$), then in $W''$ there are $6n$
connected subgraphs $X_1\l X_{6n}$, pairwise disjoint, such that
\begin{itemize}
\item for $1\le i<j\le 6n$ there is an edge of $W''$ between $X_i$ and $X_j$
\item there is no partition $(A,B,C)$ of $V(W'')$ such that 
\begin{itemize}
\item $|C|<6n$
\item no vertex in $A$ has a neighbour in $B$
\item $A$ contains at least $6n$ diagonal vertices of $W$, and
\item $B$ contains one of $X_1\l X_{6n}$.
\end{itemize}
\end{itemize}
But then, from theorem 5.4 of~\cite{GM13}, it follows that for every edge $e$ of $J_g$, with ends $j_h, j_i$ say in $J_g$,
there is a path $Q_e$ of $W''$ satisfying the following:
\begin{itemize}
\item let $e = e_{h,p} = e_{i,q}$; then the ends of $Q_e$ are $x_{h,p}$ and $x_{i,q}$;
\item the paths $Q_e\;(e\in E(J_g))$ are pairwise vertex-disjoint.
\end{itemize}
For $1\le i\le n$ let $\eta(j_i) = s_i$; and for each edge $e$ of $J_g$ with ends $j_h, j_i$ say in $J_g$, let $e = e_{h,p} = e_{i,q}$, and let $\eta(e)$
be the union of the three paths $B_{h,p},Q_e,B_{j,q}$. Thus $\eta$ is an immersion of $J_g$ in $W''$, which is itself immersed in $G$.
Consequently there is an immersion of $J_g$ in $G$ as desired.
This proves \ref{shortjumps}.~\bbox

\pagebreak

\noindent{\bf Proof of \ref{fins}.\ \ }
Let $a_1,b_1$ satisfy \ref{longjumps}, let $a_2,b_2$ satisfy \ref{inblob} and let $a_3,b_3$ satisfy \ref{shortjumps}, taking $c = a_2$.
Let $a = \max(a_1,a_2,a_3)$ and $b = 2b_1(b_2+b_3)a$.
Now let $(W,(s_i,F_i,t_i)_{1\le i\le b})$
be a fin system in $G$.
Since $s_1\l s_b$ are all diagonal vertices, we can choose $b/a=2b_1(b_2+b_3)$ of them pairwise with distance at least $a$, say $s_1\l s_{n_1}$ where $n_1 = 2b_1(b_2+b_3)$.
If some $F_{j}$ is such that at least $2b_1$ other $F_{i}$'s contain an internal vertex
of $F_j$, then the result follows from \ref{outblob}. Thus we may assume that there is no such $F_j$, and so we may assume that $F_1\l F_{n_2}$ pairwise have no internal vertex
in common, where $n_2\ge n_1/(2b_1) = b_2+b_3$.
If there are at least $b_2$ values of $i\in \{1\l n_2\}$ such that the distance from $s_i$ to $t_i$ is at least $a_2$, the result follows from \ref{inblob}, so
we assume there do not exist $b_2$ such values. But then there
are at least $b_3$ values of $i\in \{1\l n_2\}$ such that the distance from $s_i$ to $t_i$ is at most $a_2$, and the result follows from \ref{shortjumps}. 
This proves \ref{fins}.~\bbox

\section{Four-edge-connectivity}

In this section we use \ref{fins} to prove \ref{main3}. We begin with a lemma.

\begin{thm}\label{immfin}
For all $g>1$ let $b$ be as in \ref{fins}. Let $G$ be a graph, and let $\eta_0$ be an immersion in $G$ of an elementary wall $W_0$, and let $S_0$ be a set
of diagonal vertices of $W_0$, such that 
\begin{itemize}
\item if $e,f$ are distinct edges of $W_0$ and some internal vertex of $\eta_0(e)$ equals some internal vertex of $\eta_0(f)$, then there exists $s\in S_0$
incident with both $e,f$
\item $|S_0|=b$
\item for each $s\in S_0$ there is a path $F$ of $G$ with distinct ends $\eta_0(s),t$, where $t$ is a vertex of $\eta_0(W_0\setminus s)$
and no edge of $F$
belongs to $E(\eta_0(W_0)))$. 
\end{itemize}
Then there is an $\eta_0(S_0)$-rooted immersion of $J_g$ in $G$.
\end{thm}
\Proof
Let $S_0 = \{s_1\l s_b\}$, say, and for $1\le i\le b$ let $F_i$ be a path from $\eta_0(s_i)$ to some vertex $t_i$ of $\eta_0(W_0\setminus s_i)$ as in
the theorem. We proceed by induction on the sum, over all pairs of distinct edges $e,f$ of $W_0$, of the number of vertices in $\eta_0(e)\cap \eta_0(f)$.
For this quantity fixed, we proceed by induction on the sum of the lengths of $F_1\l F_b$. We may therefore
assume that for $1\le i\le b$, $t_i$ is the only vertex of $F_i$ that belongs to $\eta_0(W_0\setminus s_i)$.
If there do not exist distinct edges $e,f$ of $W_0$ such that some internal vertex
of $\eta_0(e)$ equals some internal vertex of $\eta_0(f)$, then $\eta_0(W_0)$ is a wall and the result follows from \ref{fins}. Thus we assume that some
$v$ is an internal vertex of $\eta_0(e_1)\l \eta_0(e_k)$ say, where $e_1\l e_k$ are distinct edges of $W_0$ and $k\ge 2$.
 From the hypothesis, every two of $e_1\l e_k$ have a common end in $S_0$, and since every edge has at most one end in $S_0$, it follows that $e_1\l e_k$
are all incident with some member of $S_0$, say $s_1$. It follows that $v\notin \eta_0(W_0\setminus s_1)$, and so $v\ne t_1$.

Let $d_1,d_2$ be the two edges of $\eta_0(e_1)$ incident with $v$, with ends $u_1,v$ and $u_2,v$ respectively. 
Let $G'$ be obtained from $G$
by deleting $d_1,d_2$ and adding a new edge $d_0 = u_1u_2$. Let $\eta_0'(w) = \eta(w)$ for all $w\in V(W_0)$, and $\eta_0'(d) = \eta(d)$ for all $d\in E(W_0)$
with $d\ne e_1$; let $\eta_0'(e_1)$ be a path joining the ends of $\eta_0(e_1)$ with edge set in $E(\eta_0(e_1)\setminus\{d_1,d_2\})\cup \{d_0\}$. Then $\eta_0'$
is an immersion of $W_0$ in $G'$. 
Moreover, for $1\le i\le b$, if $t_i = v$ then $e_2$ is not incident with $s_i$ (because $e_2$ is incident with none of $s_2\l s_{b_2}$, and $t_1\ne v$).
Consequently, $t_i\in \eta_0'(W_0\setminus s_i)$ for $1\le i\le b$. Thus from the inductive hypothesis, the theorem holds for $G'$, and hence it also holds for $G$.
This proves \ref{immfin}.~\bbox

\pagebreak

\noindent{\bf Proof of \ref{main3}.\ \ }
Let $b$ be as in \ref{fins} (we may assume that $b\ge 2$).
Let $W$ be a wall in $G$, and let $S$ be a set of diagonal vertices with $|S|= b$, pairwise four-edge-connected in $G$. Let $W_0$ be an elementary wall of the same height,
and let $S_0$ be the set of diagonal vertices of $W_0$ mapped to $S$ under the corresponding subdivision map. Choose an immersion $\eta_0$ of $W_0$ in $G$ and a subset
$D$ of $S_0$, with the following properties:
\begin{itemize}
\item if $e,f$ are distinct edges of $W_0$ and some internal vertex of $\eta_0(e)$ equals some internal vertex of $\eta_0(f)$, then there exists $s\in S_0$
incident with both $e,f$
\item $\eta_0(s)\in S$ for each $s\in S_0$
\item for each $s\in D$ there is a path $F$ of $G$ with distinct ends $\eta_0(s),t$, where $t$ is a vertex of $\eta_0(W_0\setminus s)$
and no edge of $F$
belongs to $E(\eta_0(W_0))$
\item subject to all these conditions, $|D|$ is maximum. 
\end{itemize}
(Satisfying the first three conditions is trivially possible, taking $D = \emptyset$ and $\eta_0$ the subdivision map to $W$.) 

Suppose that there exists $s'\in S_0\setminus D$. 
Let the three neighbours of $s'$ in $W_0$ be $x_1,x_2,x_3$, and let $B_1,B_2,B_3$ be the images under $\eta_0$ of
the edges $s'x_1,s'x_2,s'x_3$. Now since $|S|\ge 2$ and $\eta_0(s')$ is 
four-edge-connected to the other members of $S$, it follows that there are four edge-disjoint paths $P_1\l P_4$
of $G$ from $\eta_0(s')$ to $V(\eta_0(W_0\setminus s'))$, each with no internal vertex in $V(\eta_0(W_0\setminus s'))$; and since the three branches
of $\eta_0(W_0)$ incident with $\eta_0(s')$ provide three such paths, it follows from the theory of augmenting paths that $P_1\l P_4$
can be chosen such that $P_1,P_2,P_3$ have final vertices $\eta_0(x_1), \eta_0(x_2),\eta_0(x_3)$ respectively. Let $P_4$ have final vertex
$t'$ say. Now let
$\eta_0'(v) =\eta_0(v)$ for each $v\in V(W_0)$, and $\eta_0'(e) = \eta_0(e)$ for every edge $e$ of $W_0$ not incident with $s'$; let $\eta_0'(s'x_i) = P_i$ for
$i = 1,2,3$. It follows that $\eta_0'$ is an immersion of $W_0$ in $G$, satisfying the second condition above. Moreover, the first condition above is satisfied;
for if $e,f$ are distinct edges of $W_0$ and some internal vertex of $\eta_0'(e)$ equals some internal vertex of $\eta_0'(f)$, we may assume that
$e$ is incident with $s'$ and $f$ is not; but then some internal vertex of one of $P_1,P_2,P_3$ belongs to $\eta_0(W_0\setminus s')$, a contradiction.
Let $D' = D\cup \{s'\}$; we claim that the third condition is satisfied. Certainly, the new member $s'$ of $D'$ satisfies the condition,
since $P_4$ is a path of $G$ with distinct ends $\eta_0'(s'),t'$, 
and $t'$ is a vertex of $\eta_0'(W_0\setminus s)$
and no edge of $P_4$
belongs to $E(\eta_0'(W_0)))$; but we must check that the other members of $D'$ still satisfy the condition. Thus, let $s\in D$.
There is a path $F$ of $G$ 
with distinct ends $\eta_0(s),t$, where $t$ is a vertex of $\eta_0(W_0\setminus s)$
and no edge of $F$ belongs to $E(\eta_0(W_0))$. The same path $F$ also works for $\eta_0'$ unless either
\begin{itemize}
\item $t$ is not a vertex of $\eta_0'(W_0\setminus s)$; this implies that $t$ is an internal vertex of one of $B_1,B_2,B_3$, and then we can augment $F$ to a path
with the desired properties; or
\item some edge or internal vertex of $F$ belongs to $E(\eta_0'(W_0)))$; but this implies that some edge or internal vertex of $F$ belongs to one of $P_1,P_2,P_3$,
and then a subpath of $F$ has the desired properties, since $s\ne s'$.
\end{itemize}
Thus the three conditions are still satisfied, contrary to the maximality of $|D|$. Thus there is no such $s'$, and so 
$D = S_0$, and the result follows from \ref{immfin}. This completes the proof of \ref{main3}.~\bbox

\pagebreak

\end{document}